\magnification=1100
\hsize=4.7 true in
\vsize=7.5 true in
\hoffset=.9 true in
\voffset=.8 true in 
\baselineskip=12.65pt
\def\p{\vskip .2cm}

\def\C{{\bf C}}
\def\K{{\bf K}}
\def\Q{{\bf Q}}
\def\F{{\bf F}}
\def\Z{{\bf Z}}

\def\N{{\bf N}}

\def\O{\Omega} 
\def\div{{\rm div\,}}
\def\Exc{{\rm Exc\,}}
\def\Spec{{\rm Spec\,}}

\centerline{\bf ON THE TOPOLOGY OF BIRATIONAL MINIMAL MODELS}\p\p\p\p

\centerline{Chin-Lung Wang}\p\p\p\p\p\p
\hfill{\it Dedicated to Hui-Wen and Chung-Yang}
\p\p\p\p\p\p\p\p

\noindent
{\bf Introduction}\p

In the study of higher dimensional algebraic geometry, an important 
reduction step is to study certain good birational models of a given 
algebraic manifold. This leads to the famous ``minimal model program'' 
initiated by Mori -- the search for birational models with numerically 
effective canonical divisors and with at most terminal singularities.
The existence problem is still completely open in dimensions bigger
than three, but even worse, in contrast to the two dimensional case, 
the minimal model is not unique in higher dimensions. It is then an
important question to see what kind of invariants are shared by
all the birationally equivalent minimal models. And more generally,
to see what kind of invariants are preserved under certain elementary 
birational transformations. In this paper, some results in this 
direction are given:\p

{\bf Theorem A.} {\it Let $f\colon X\cdot\!\cdot\!\!\to X'$ be a
birational map between two smooth complex projective varieties such
that the canonical bundles are numerically effective along the 
exceptional loci, then $X$ and $X'$ have the same Betti numbers.
In particular, birational smooth minimal models have the same
Betti numbers.}\p

Theorem A generalizes, in the smooth case, previous results of Koll\'ar 
on the invariance of cohomologies under flops in dimension three
(cf.\ 5.1). Another interesting corollary via the Mayer-Vietoris 
argument shows that the exceptional loci of the given birational map 
also share the same Betti numbers (Corollary 4.5).\par
 
The proof of Theorem A is based on general considerations in birational 
geometry given in \S1 and Grothendieck-Deligne's solution to the 
Weil conjecture [D1, D2]. The bridge to connect these two is the 
theory of $p$-adic integrals.\par

The idea to use the Weil conjecture via
$p$-adic integrals to compute cohomologies can be dated back to 
Harder and Narasimhan in the 70's [HN]. However, it was used there in a
somewhat different way. Recently this was taken up again by Batyrev by
developing Weil's idea of $p$-adic measure [Ba]. In fact, he established 
Theorem A in the special case of projective Calabi-Yau manifolds. \par

By extending this idea further, an argument based on birational
correspondences is developed here in order to deal with the general 
case. Namely, we introduce in \S1 the notion of ``K-partial ordering'' 
and relate it to interesting geometric situations. The applicability of the
Weil conjecture is largely clarified in terms of this notion (cf.\
Proposition 2.16 and Theorem 3.1). Moreover, this approach also provides 
a natural setting in the singular case. \par

In this paper, We have tried to develop this, together with the $p$-adic 
measure, as far as possible so that it could fit the need of the minimal 
model theory. In fact, an easy but very interesting result observed here
is that the integral points of a $p$-adic variety has finite $p$-adic
measure if and only if it has at most log-terminal singularities 
(Proposition 2.12). This gives the basic reason why $p$-adic integrals
fit into the framework of minimal model theory naturally.  
But due to technical reasons, we have restricted ourself 
to the smooth case when we state and prove Theorem A. 
(See however 5.3 for the singular case.)\par

Theorem A is still not all satisfactory in two aspects -- 
the torsion elements are not considered and no natural maps between 
cohomologies has been mentioned. Although there is one obvious candidate 
for this map -- the cohomology correspondence induced from the birational 
correspondence, it is not clear how to show directly
that it induces isomorphisms. In fact, there is no strong evidence why
this should be true. The next result only deals with the simplest 
cases. However, it is included to emphasize this important aspect.\p 

{\bf Theorem B.} {\it Smooth minimal models minimize $H^2(X,\Z)$ 
compatible with the Hodge structure among birational smooth projective
varieties. In the singular case, the minimal models minimize the group 
of Weil divisors among birational projective varieties with 
at most terminal singularities.}\p

The proof, which is elementary (does not use the Weil conjecture), 
is contained in $\S4$ together with some related results. In fact, it
is simply another application of the notion of K-partial ordering.\par 

To finish the introduction, it is worth pointing out 
that in stating both theorems, what we have in mind is that 
there should be a ``minimal cohomology theory'' among birational varieties.
Moreover, it should be realized exactly by the minimal models.
\p\p\p

\noindent
{\bf\S1 Birational Geometry}\p

We begin with some standard definitions. For a complete treatment
of minimal model theory, the reader should consult [KMM].\par

Let $X$ be an $n$ dimensional complex normal $\Q$-Gorenstein variety. 
That is, the canonical divisor $K_X$ is $\Q$-Cartier. Recall that $X$ has 
(at most) terminal (resp.\ canonical, resp.\ log-terminal) singularities
if there is a resolution $\phi:Y\to X$ such that in the canonical 
bundle relation 
$$
K_Y=_\Q \phi^*K_X+\sum a_iE_i, \leqno (1.1)
$$ 
we have that $a_i>0$ (resp.\ $a_i\ge 0$, resp.\ $a_i>-1$) for all $i$.
Here, the $E_i$'s vary among the prime components 
of all the exceptional divisors. Although (1.1) holds only up to \Q-linear
equivalence, the divisor $\sum a_i E_i\in Z_{n-1}\otimes\Q$ is uniquely
determined. Moreover, the condition on $a_i$'s is readily seen to be
independent of the chosen resolution. It is also elementary to see that 
smooth points are all terminal.\par

Let $Z$ be a proper subvariety of $X$. A $\Q$-Cartier divisor $D$ is
called numerically effective (nef) along $Z$ if 
$D.C:=\deg_{\tilde C}(f^*D)\ge 0$ for all effective curves $C\subset Z$, 
where $f\colon \tilde C\to C$ is the normalization of $C$. $D$ is simply 
called nef if $Z=X$. A projective variety $X$ is called a minimal model 
if $X$ is terminal and $K_X$ is nef.\par

Two normal varieties $X$ and $X'$ are birational if they have
isomorphic function fields $K(X)\cong K(X')$ (over \C). Geometrically, this
means that there is a rational map $f\colon X\cdot\!\cdot\!\!\to X'$ 
such that $f^{-1}$ is also rational. The exceptional loci of $f$ are 
defined to be the smallest subvarieties $Z\subset X$ and $Z'\subset X'$
such that $f$ induces an isomorphism $X-Z\cong X'-Z'$.\par

Among the class of birational \Q-Gorenstein varieties, We have the notion of 
{\bf K-partial ordering} (where the ``K'' is for canonical divisors):\p

{\bf Definition 1.2.}  For two \Q-Gorenstein varieties $X$ and $X'$, 
we say that $X\le_K X'$ (resp.\ $X <_K X'$) if there is a birational
correspondence $(\phi,\phi'): X\leftarrow Y \to X'$ with $Y$ smooth, 
such that $\phi^*K_X \le_\Q \phi'^*K_{X'}$ (resp.\ ``$<_\Q$''). 
Moreover, ``$X\le_K X'$'' plus ``$X\ge_K X'$'' implies that ``$X=_K X'$'',
ie. $\phi^*K_X =_\Q \phi'^*K_{X'}$. In this case, we say that 
$X$ and $X'$ are K-equivalent.\p

The well-definedness of this notion follows from the canonical bundle 
relations
$$
K_Y =_\Q \phi^*K_X + E =_\Q \phi'^*K_{X'} + E', \leqno (1.3)
$$
since we know that $X\le_K X'$ if and only if $E\ge E'$. In the terminal case, 
this means that $\phi$ has more exceptional divisors than $\phi'$ 
(so heuristically, $X$ is ``smaller'' than $X'$).\par 

Here is the typical geometric situation that we can compare their
K-partial order:\p

{\bf Theorem 1.4.} {\it Let $f\colon X\cdot\!\cdot\!\!\to X'$ 
be a birational map between two varieties with canonical singularities.
Suppose that the exceptional locus $Z\subset X$ is proper and that 
$K_X$ is nef along $Z$, then $X\le_K X'$. Moreover, if $X'$ is
terminal, then $Z$ has codimension at least two.}\p

{\it Proof.} Let $\phi:Y\to X$ and $\phi':Y\to X'$
be a good common resolution of singularities of $f$ so that 
the union of the exceptional set of $\phi$ and $\phi'$ is a
normal crossing divisor of $Y$. This can be done by 
considering $\bar{\Gamma}_f\subset X\times X'$,
the closure of the graph of $f$, blowing up the exceptional set of
$\bar{\Gamma}_f\to X$ and $\bar{\Gamma}_f\to X'$ and then
taking $Y$ to be a Hironaka (embedded) resolution [Hi].\par

Consider the canonical bundle relations:
$$
\eqalign{K_Y&=_\Q \phi^*K_X + E \equiv \phi^*K_X + F + G\cr
	    &=_\Q \phi'^*K_{X'} + E' \equiv \phi'^*K_{X'} + F' + G'.\cr} 
\leqno (1.5)
$$
Here $F$ and $F'$ denote the sum of divisors (with coefficients $\ge 0$) 
which are both $\phi$ and $\phi'$ exceptional. $G$ (resp.\ $G'$) denotes 
the part which is $\phi$ exceptional but not $\phi'$ exceptional
(resp.\ $\phi'$ but not $\phi$ exceptional). Notice that 
$\phi(G')\subset Z$.\par

To proceed, we write
$$
\phi'^*K_{X'} =_\Q \phi^* K_X + G + (F - F' - G').\leqno (1.6)
$$
It is enough to prove that $F-F'-G'\ge 0$, because this implies 
that $F-F'\ge 0$ and $G'=0$, and so $E\ge E'$.\par

By taking a generic hyperplane section $H$ of $Y$ $n-2$ times,
the problem is reduced to a problem on surfaces. Namely
$$
H^{n-2}.\phi'^*K_{X'}=_\Q H^{n-2}.\phi^*K_X 
  + \zeta + (\xi-\xi'-\zeta'), \leqno (1.7)
$$
where $\xi = H^{n-2}.F$ and $\zeta = H^{n-2}.G$ etc. If $\xi-\xi'-\zeta'$ 
is not effective, write it as $H^{n-2}.(A-B)=a-b$ with $A$ and $B$ 
effective. Then by taking the intersection of (1.5) with $b$, we get
$$
B.H^{n-2}.\phi'^*K_{X'}=_\Q B.H^{n-2}.\phi^*K_X 
  + b.\zeta + b.a -b^2.\leqno (1.8)
$$
The left hand side is zero since $B\subset E'$ is $\phi'$ exceptional. 
Moreover, if $B\subset F'$ then $B.H^{n-2}.\phi^*K_X=0$ too. If
$B\subset G'$ then the curve $\phi(B.H^{n-2})\subset \phi(G')
\subset Z$ is inside the exceptional locus. So the first three
terms in the right hand side are non-negative since $K_X$ is nef 
along $Z$ and $a$, $b$ and $\xi$ are different components. However, 
since $b$ is a nontrivial combination of exceptional curves in 
$H^{n-2}$, we have from the Hodge index theorem that $b^2<0$, a 
contradiction! Hence $F-F'-G'\ge 0$.\par

For the second statement, from the construction of $Y$, we know that all 
components of the exceptional sets, denoted by $\Exc\phi$ and 
$\Exc\phi'$ respectively, are divisors. If $X'$ is assumed to be 
terminal, then all $\phi'$ exceptional divisors occur as 
components of $E'$. So $G'=0$ implies that $\Exc\phi'\subset\Exc\phi$.
With this understood, from 
$$
X-\phi(\Exc\phi)\cong Y-\Exc\phi\cong X'-\phi'(\Exc\phi)\subset 
X'-\phi'(\Exc\phi'),\leqno (1.9)
$$
we conclude that $Z\subset\phi(\Exc\phi)$ is of codimension at least two. 
Q.E.D.\p

{\bf Corollary 1.10.} {\it Let $f\colon X\cdot\!\cdot\!\!\to X'$ 
be a birational map between two terminal varieties such that 
$K_X$ (resp. $K_{X'}$) is nef along the exceptional locus $Z\subset X$
(resp. $Z'\subset X'$), then $X=_K X'$ and $f$ is an isomorphism in 
codimension one. This applies, in particular, if both $X$ and $X'$ are 
minimal models.}\p 

{\bf Variant 1.11.} Instead of assuming that the exceptional locus in
$X$ is proper, one can generalize Theorem 1.4 to the relative case, namely 
$f$ is a $S$-birational map and that $X\to S$ and $X'\to S$ are proper 
$S$-schemes. The proof is identical to the one given above by 
changing notation.\p

{\bf Remark 1.12.} This type of argument is familiar in the minimal
model theory. Notably, in analyzing the log-flip diagram 
(eg. [KMM; 5-1-11]) or more specially, the flops. Theorem 1.4 implies that
if $X'$ is a flip of $X$, then $X \ge_K X'$ (in fact, more is true:
$X >_K X'$). Corollary 1.10 implies that flop induces K-equivalence. 
Since flip/flop will not be used in any essential way in this paper, 
we will refer the interested reader to [KMM] for the definitions.\par

The proof given above is inspired by Koll\'ar's treatment of flops in [Ko].
\p\p\p

\noindent
{\bf \S2 The Weil Conjecture and $p$-adic Integrals}\p

To prove Theorem A, we will show that $X$ and $X'$ 
have the same number of rational points over certain finite fields 
when a suitable good reduction is taken. That is, we prove that they 
have the same ``zeta function''. The theorem will then follow from the
statement of the Weil conjecture. \p

{\bf 2.1. The reduction procedure.} This is 
standard in algebraic geometry and in number theory: as long as we 
perform only a finite number of ``algebraic constructions''
in the complex case, e.g. consider morphisms, since all the 
objects involved can by defined by a finite number of polynomials, 
we can take $S\subset\C$ a finitely generated subring over $\Z$ so 
that everything is defined over $S$. $S$ has the property that
the residue field $S/m$ of any maximal ideal $m\subset S$ is finite.\par

If we start with ``smooth objects'', general reduction theory then says 
that for an infinite number of 
``good primes'' (in fact, Zariski dense in $\Spec (S)$), we may 
get good reductions so that everything is defined 
smoothly over the finite residue field $\F_q$ with $q=p^r$ for some prime 
number $p$. We may also assume that this reduction has a lifting such
that everything is defined smoothly over $R$, the maximal compact 
subring of a $p$-adic local field $\K$, i.e. a finite extension
field of $\Q_p$, with residue field $\F_q$.\par

Here is a special way to see this. Let $F$ be the quotient field of $S$.
Based on the fact (and others) that $\Q_p$ has infinite transcendence 
degree, the ``embedding theorem'' (see for example 
[Ca; p.82]) says that for an infinite number of $p$'s, there is an 
embedding of fields $i:F\to \Q_p$ such that $i(S)\subset \Z_p$. Moreover,
$i$ may be chosen so that a prescribed finite subset of $S$, say 
the coefficients of those defining polynomials, is mapped into the 
set of $p$-adic units. This embedding then gives the desired lifting.\par 

Let $P$ be the unique maximal ideal of $R$ (so $R/P\cong \F_q$). 
We denote by $\bar X$, $\bar U, \ldots$ those objects constructed 
from $X$, $U \ldots$ via reductions mod $P$. That is, objects 
lie over the point $\Spec R/P\to \Spec R$ -- they are defined over $\F_q$. 
We also denote the reduction map by $\pi:X(R)\to \bar X(\F_q)$ etc.\p  

{\bf 2.2. The Weil conjecture.}
Let $\bar X$ be a variety defined over a finite field $\F_q$.
After fixing an algebraic closure, the Weil zeta function of $\bar X$ 
is defined by
$$
Z(\bar X,t):=\exp\Bigg(\sum\nolimits_{k\ge 1}|\bar X(\F_{q^k})|
{t^k\over k}\Bigg).\leqno (2.3)
$$\par

In 1949, Weil conjectured several nice properties of this zeta 
function for smooth projective varieties and expalined how some of 
these would follow once a suitable cohomology theory exists [W1].
This lead Grothendieck to his creation of \'etale cohomology theory.\par

More precisely, Grothendieck proved a ``Lefschetz fixed point formula'' 
in a very general context (eg. constructible sheaves over seperated
schems of finite type) [D2], which in particular
implies that the zeta function is a rational function:
$$
Z(\bar X,t)={P_1(t)\cdots P_{2n-1}(t)
\over P_0(t)P_2(t)\cdots P_{2n}(t)},\leqno (2.4)
$$
where $P_j(t)$ is a polynomial with integer coefficients 
such that $P_j(0)=1$ and ${\rm deg}\,P_j(t)=h^j$, the $j$-th Betti number
of the compactly supported $\ell$-adic \'etale cohomologies (for a prime
number $\ell\ne p$). Moreover, when $\bar X$ comes from a good 
reduction of a smooth complex projective variety $X$ in the sense described
in (2.1), $h^j$ coincides with the $j$-th Betti number of the singular 
cohomologies of $X(\C)$.\par 

Deligne [D1] completed the proof of the Weil conjecture by proving the 
important ``Riemann Hypothesis'' that all roots of $P_j(t)$ have 
absolute value $q^{-j/2}$. In particular, the complete information 
about the $\F_{q^k}$-rational points determines 
the $h^j$'s and all the roots.\p

{\bf 2.5. Counting points via $p$-adic integrals.}
How do we count $\bar X(\F_q)$? If $\bar X$ comes from the good 
reduction of a smooth 
$R$-scheme, we will see that such a counting can be achieved by 
using $p$-adic integrals (cf.\ Theorem 2.8). We will first 
recall some elementary aspects of the $p$-adic integral over 
$\K$-analytic manifolds and over $R$-schemes. \par
 
Consider the Haar measure on the locally compact field $\K$,
normalized so that the compact open ``disk'' $R$ has volume 1:
$$
\int\nolimits_R |dz|=1. \leqno (2.6)
$$
We may extend this to the multivariable case and define the $p$-adic 
integral of any regular $n$ form $\Psi=\psi(z_1,\cdots,z_n)dz_1\wedge
\cdots\wedge dz_n$ by
$$
\int\nolimits_{R^n} |\Psi|:=\int\nolimits_{R^n}|\psi(z)| 
|dz_1\wedge\cdots\wedge dz_n|.\leqno (2.7)
$$
Here $|a|:=q^{-\nu_p({\rm N}_{\K/\Q_p}(a))}$ is the usual $p$-adic norm.\par

We may define an integral slightly more general than (2.7):
suppose that $\Psi$ is a $r$-pluricanonical form such that in local
analytic coordinates $\Psi=\psi(z_1,\cdots,z_n)(dz_1\wedge
\cdots\wedge dz_n)^{\otimes r}$. We define the integration of a
``$r$-th root of $|\Psi|$'' by
$$
\int\nolimits_{R^n} |\Psi|^{1/r}:=\int\nolimits_{R^n}
|\psi(z)|^{1/r} |dz_1\wedge\cdots\wedge dz_n|.\leqno (2.7')
$$
This is independent of the choice of coordinates, as can be checked
easily by the same method as in [W2; p.14]. So we can extend the 
definition to (not necessarily complete) $\K$-analytic manifolds 
with $\Psi$ a (possibly meromorphic) pluricanonical form.
Certainly then the integral defined may not be finite.\par

The key property we need is the following (slightly more general form 
of a) formula of Weil [W2; 2.2.5]. We briefly sketch its proof.\p 

{\bf Theorem 2.8.} {\it Let $U$ be a smooth $R$-scheme and $\O$ a nowhere
zero $r$-pluricanonical form on $U$, then
$$
\int\nolimits_{U(R)}|\O|^{1/r}={|\bar U(\F_q)|\over q^n}.
$$}\p

{\it Proof.} The proof given by Weil in [W2] goes through
without difficulties -- one first observes that the
reduction map $\pi\colon U(R)\to \bar U(\F_q)$ induces an isomorphism 
between $\pi^{-1}(\bar t)$ and $PR^n$ for any $\bar t\in \bar U(\F_q)$ 
(Hensel's lemma) such that there is
a function $\psi$ with $|\psi(z)|=1$ and 
$$
\O=\psi(z)\cdot(dz_1\wedge\cdots\wedge dz_n)^{\otimes r} 
\leqno (2.9)
$$ 
in the $\K$-analytic chart $PR^n$. This implies that
$\int_{\pi^{-1}(\bar t)}|\O|^{1/r}=1/q^n$
for any $\bar t \in \bar U(\F_q)$. Summing over $\bar t$ then gives the
result. Q.E.D.\p

The right hand side of (2.8) shows that the integral is independent of
the choice of the form $\O$. One may also see this by observing that 
any two such forms differ by a nowhere vanishing function on $U$ 
(over $R$) which takes values in the units on all $R$-points. 
This allows one to define a canonical $p$-adic measure on the 
$R$-points of smooth $R$-schemes by ``gluing'' the local integrals.
We will define it in the singular case with the hope that it may be useful 
for later development.\p

{\bf 2.10. Canonical measure on $\Q$-Gorenstein $R$-schemes.}
We will only consider those $R$-schemes, eg. $X$, that come from complex 
\Q-Gorenstein varieties as in (2.1). Let $r\in \N$ such that $rK_X$ 
is Cartier (locally free). We may assume that we have a $R$-resolution
of singularities $\phi\colon Y\to X$, which is a projective $R$-morphism, 
so that the reduced part of the exceptional
set is a simple normal crossing $R$-variety. We will define a measure 
on $X(R)$ such that the measurable sets are exactly the compact open 
subsets in the \K-analytic topology.\par

Let $U_i$'s be a Zariski open cover of $X$ such that $rK_X|_{U_i}$ is 
actually free. Then for a compact open subset 
$S\subset U_i(R)\subset X(R)$, we define its measure by
$$
m_X(S)\equiv \int\nolimits_S |\O_i|^{1/r}
:= \int\nolimits_{\phi^{-1}(S)} |\phi^*\O_i|^{1/r}, \leqno (2.11)
$$
where $\O_i$ is an arbitrary generator of $rK_X|_{U_i}$. Notice that 
the properness of $\phi$ implies that $\phi^{-1}(S)\subset Y(R)$. This
allows us to operate the integral entirely on $R$-points.\par

For general compact open $S\subset X(R)$, we may break $S$ into 
disjoint pieces 
$S_j$ so that $S_j$ is contained in some $U_i(R)$ (in fact, $S_j$ may be 
chosen to lie entirely in a fiber of the reduction map $\pi$), and then
define $m_X(S)=\sum_i m_X(S_i)$. Notice that $m_X(S)$ is again independent 
of the choice of $U_i$, $\O_i$ and $Y$.\par

The following proposition explains the possible connection between
the canonical measure and the minimal model theory:\p 

{\bf Proposition 2.12.} {\it For a \Q-Gorenstein $R$-variety $X$, 
$X(R)$ has finite measure if and only if $X$ has at 
most log-terminal singularities.}\p

{\it Proof.} Consider the canonical bundle relation for $\phi\colon Y\to X$
$$
rK_Y = \phi^* rK_X + \sum\nolimits_i e_i E_i \leqno (2.13)
$$
with $rK_X$ being Cartier and $e_i\in\Z$. To determine the finiteness 
of $m_X(X(R))$, we only
need to consider those $R$-points on the exceptional fibers. Locally, 
$\div\phi^*\O=\sum_i e_i E_i$ for a generator $\O$ of $rK_X$. So
the integral is a product of one dimensional integrals of the form 
$$
I_i := \int\nolimits_R |z^{e_i}\,dz^{\otimes r}|^{1/r}
= \int\nolimits_R |z|^{e_i/r}\,|dz|. \leqno (2.14)
$$
If this is finite, then
$$
I_i = \int\nolimits_{PR} |z|^{e_i/r}\,|dz| + (q - 1){1 \over q}
= q^{-(e_i/r + 1)}I_i + {q - 1 \over q}. \leqno (2.15)
$$
Since $I_i > 0$, this makes sense only if $q^{e_i/r + 1} > 1$. 
That is, $e_i/r > -1$, which is exactly the definition of 
log-terminal singularities. Q.E.D.\p

Since the measure is defined Zariski-locally via $p$-adic integrals, 
for smooth $X$, we have from Weil's formula (2.8) that:\p

{\bf Corollary 2.16.} {\it Let $X$ be an $n$-dimensional smooth 
$R$-variety with finite residue field $\F_q$, then
$$
m_X(X(R)) = {|\bar X(\F_q)| \over q^n}.
$$}\p

{\bf Remark 2.17.} If $X$ is singular, $m_X((X(R))$ is a weighted 
counting of the rational points. By definition, if $\phi\colon Y\to X$ 
is a crepant $R$-morphism, ie. $K_Y=_\Q\phi^*K_X$, then 
$m_X((X(R))=m_Y((Y(R))$. In particular, $m_X((X(R))$ counts the rational 
points of $\bar Y$ if $Y$ is smooth! This applies to many interesting 
``pure canonical'' singularities and to terminal singularities having 
small resolutions. However, further investigation on the 
precise ``geometric meaning'' of this weighted counting is still needed 
for the general case (cf.\ 5.3).
\p\p\p

\noindent
{\bf \S3 The Proof of Theorem A}\p

We will in fact prove a result which connects the notion of K-partial
ordering and the canonical measure. This will largely clarify the role
played by the Weil conjecture.\p

{\bf Theorem 3.1.} {\it Let $X$ and $X'$ be two birational log-terminal 
$R$-varieties. Then $m_X(X(R)) \le m_{X'}(X'(R))$ if $X \le_K X'$.}\p

{\it Proof.} Consider as before, a birational correspondence
$(\phi,\phi'): X\leftarrow Y \to X'$ over $R$ with $Y$ a smooth 
$R$-variety. Let $r \in \N$ be such that both $rK_X$ and $rK_{X'}$ are
Cartier. Then $X \le_K X'$ if and only if 
in the canonical bundle relations
$rK_Y = \phi^*rK_X + E = \phi'^*rK_{X'} + E'$, we have $E \ge E'$.\par

From the properness of $\phi$ and $\phi'$, 
we have that $\phi^{-1}(X(R))=Y(R)=\phi'^{-1}(X'(R))$. So from the 
definition of the measure (2.11), it suffices to show that for any 
compact open subset $T\subset Y(R)$ with $\pi(T)$ a single point 
$\bar y \in \bar Y (\F_q)$, we have
$$
\int\nolimits_{T} |\phi^*\O|^{1/r}
\le \int\nolimits_{T} |\phi'^*\O'|^{1/r}. \leqno (3.2)
$$
Here $\O$ is an arbitrary local generator of $rK_X$ on a Zariski open
set $U$ where $rK_X$ is actually free and such that 
$\bar\phi(\bar y)\in \bar U$ (and with similar conditions for $\O'$).\par
 
Clearly, (3.2) can fail to be an equality only if 
$\bar y\in \bar E\cup \bar E'$. However, in this case $E\ge E'$ says that
the order of $\phi^*\O$ is no less than that of $\phi^*\O$.
(3.2) then follows from the definition of the $p$-adic integral $(2.7')$
(see also (2.15)). Q.E.D.\p

If $X$ and $X'$ are smooth, combining this with (2.16) gives \p

{\bf Corollary 3.3.} {\it Let $X$ and $X'$ be two birational smooth 
$R$-schemes. Then $|\bar X(\F_q)| \le |\bar X'(\F_q)|$ if $X \le_K X'$.}\p

With this been done, by working on cyclotomic extensions of $\K$, 
the same proof shows that $|\bar X(\F_{q^k})| \le |\bar X'(\F_{q^k})|$ 
for all $k \in \N$. In particular, $Z(\bar X, t)\le Z(\bar X', t)$
for all $t>0$. The same is true for all the derivatives, but
it is not clear how to make use of these.\p

{\bf Corollary 3.4.} {\it Let $X$ and $X'$ be two birational complex
smooth varieties. They have the same Euler number for the compactly
supported cohomologies if $X =_K X'$.}\p

{\it Proof.} Apply the reduction procedure (2.1) to reduce this 
to the $p$-adic case. The statement then follows from Grothendieck's 
Lefschetz fixed point formula (2.4) and the above comparison of zeta 
functions. Q.E.D.\p

What kind of geometric situation can we have $X\le_K X'$? Theorem 1.4
provides such a typical case inspired by the minimal model theory.
Namely, let $f\colon X\cdot\!\cdot\!\!\to X'$ be a birational map between 
two varieties with at most canonical singularities and with proper 
exceptional locus $Z\subset X$ such that $K_X$ is nef along $Z$, 
then $X\le_K X'$.\par

So far we have not used Deligne's theorem on the ``Riemann Hypothesis''.
To use it, we need to impose the projective assumption.\p
  
{\bf Theorem 3.5.} {\it Let $X$ and $X'$ be two birational smooth 
projective $R$-schemes. If $X =_K X'$ then $m_X(X(R)) = m_{X'}(X'(R))$.
Equivalently, $Z(\bar X, t) = Z(\bar X', t)$.
In particular, they have the same \'etale Betti numbers
by the Weil conjecture.}\p

Now Theorem A simply follows from the reduction procedure (2.1),
Corollary 1.10 and Theorem 3.5. Q.E.D.\p

{\bf Remark 3.6.} In the preliminary version of this paper (dated
October 1997), Theorem A was stated with the assumption that
the canonical bundle is semi-ample, that is, $rK_X$ is generated by
global sections for some $r\in \N$. The proof proceeds by cutting
out the pluri-canonical divisors and applying $p$-adic integrals
to the birational correspondence, where the notion of K-equivalence
is essential for this step to work.\par

By using Weil's formula (2.8),
the proof is then concluded by induction on dimensions. In this 
approach, the usage of integration of a $r$-th root of the absolute 
value of a pluricanonical form was suggested to the author by C.-L. Chai 
in order to deal with the case that $r > 1$.
Happily enough, as the author realized later, the semi-ample 
assumption can be removed once we observed that the problem can even 
be localized to the exceptional loci.\p

{\bf Remark 3.7.} The equivalence of 
zeta functions is a stronger statement than the equivalence of 
Betti numbers. Moreover, we have in fact established the equivalence 
of zeta functions for a dense set of primes. From the theory of motives, 
this suggests that we may in fact have the equivalence of Hodge structures. 
Further investigation in this should be interesting and important.\p

{\bf Question 3.8.} Is Theorem A true for K\"ahler manifolds?
\p\p\p

\noindent
{\bf \S4 Miscellaneous Results and The Proof of Theorem B}\p

Now we come back to the complex number field and begin with an 
elementary observation:\p

{\bf Lemma 4.1.} {\it If the exceptional loci of a birational map 
$f\colon X\cdot\!\cdot\!\!\to X'$ between two smooth projective varieties 
have codimension at least two then for $i\le 2$ we have 
$\pi_i(X)\cong\pi_i(X')$ and $H^i(X,\Z)\cong H^i(X',\Z)$ which 
is compatible with the rational Hodge structures.}\p
 
{\it Proof.} The real codimension four statement plus the transversality
argument shows that $\pi_i(X)\cong\pi_i(X')$, $H_i(X,\Z)\cong H_i(X',\Z)$ 
and $H^i(X,\Z)\cong H^i(X',\Z)$ canonically for $i\le 2$. 
Moreover, by Hartog's extension we know that the Hodge groups
$H^0(\O^i)$ are all birational invariants among smooth varieties.
The orthogonality of Hodge filtrations then shows 
that $H^i(X,\Q)$ and $H^i(X',\Q)$ share the same 
rational Hodge structures for $i\le 2$. Q.E.D. \p

A slightly deeper result is given by \p

{\bf Proposition 4.2.} {\it If the exceptional loci $Z\subset X$ and
$Z'\subset X'$ of a birational map $f$ between two smooth varieties
 have codimension at least two, 
then $h^i(X) - h^i(Z) = h^i(X') - h^i(Z')$.}\p

{\it Proof.} Construct a birational correpondence $X \leftarrow Y \to X$ 
as in \S1 and denote the exceptional divisor of $\phi\colon Y\to X$ (resp.\ 
$\phi'\colon Y\to X'$) by $E$ (resp.\ $E'$). Since Hironaka's resolution 
process only blows up smooth centers inside the singular set of the graph 
of $f$, the isomorphism $X-Z\cong X'-Z'$ implies that 
$\phi(E\cup E')\subset Z$ and $\phi'(E\cup E')\subset Z'$, hence
that $E_{\rm red}=E'_{\rm red}$, $Z=\phi(E)$ and $Z'=\phi'(E')$.\par
  
Consider an open cover $\{V, W\}$ of $X$ by letting $V:=X-Z$ and 
$W\supset Z$ be a deformation retract neighborhood. Let 
$\tilde V:=\phi^{-1}(V)$ and $\tilde W:=\phi^{-1}(W)\supset E$ be 
the corresponding open cover of $Y$. Then we have the following 
commutative diagram of integral cohomologies
$$
\matrix{
H^{i-1}(\tilde V\cap\tilde W) &\to &H^i(Y) &\to
  &H^i(\tilde V)\oplus H^i(E) &\to &H^i(\tilde V\cap\tilde W) \cr
\uparrow & &\uparrow & &\uparrow & &\uparrow \cr
H^{i-1}(V \cap W) &\to &H^i(X) &\to
  &H^i(V)\oplus H^i(Z) &\to &H^i(V \cap W) \cr}
\leqno (4.3)
$$\par

It is a general fact that $\phi^*\colon H^i(X)\to H^i(Y)$ is injective
(by the projection formula, that $\phi$ is proper of degree one implies 
that $\phi_!\circ\phi^*(a)=a$ for all $a\in H^i(X)$). Since 
$\tilde V \cong V$ and $\tilde V\cap\tilde W \cong V \cap W$,
simple diagram chasing shows that $H^i(Z) \to H^i(E)$ is also 
injective. We may then break (4.3) into short exact sequences 
$$
0 \to \phi^* H^i(X) \to H^i(Y) \to 
  H^i(E)/\phi^* H^i(Z) \to 0.
  \leqno (4.4)
$$
Similarly, we have for $\phi'\colon Y\to X'$:
$$
0 \to \phi'^* H^i(X') \to H^i(Y) \to 
  H^i(E')/\phi'^* H^i(Z') \to 0.
  \leqno (4.4')
$$
Since $E_{\rm red}=E'_{\rm red}$, the proposition follows immedeately.
Q.E.D. \p

Combining this with Theorem A gives \p

{\bf Corollary 4.5.} {\it Let $f\colon X\cdot\!\cdot\!\!\to X'$ be a
birational map between two smooth complex projective varieties such
that the canonical bundles are numerically effective along the 
exceptional loci, then the exceptional loci also have the same Betti 
numbers. In particular, this applies to birational smooth minimal 
models.}\p

{\bf Remark 4.6.} The proof of Theorem A in fact also shows that $\bar Z$ 
and $\bar Z'$ have the same number of $\F_q$-rational points. This is 
simply because $|\bar X(\F_q)|=|\bar X'(\F_q)|$ and 
$\bar X-\bar Z\cong\bar X'-\bar Z'$. In particular, if $Z$ and $Z'$ are 
smooth then they have the same Betti numbers. Although this argument 
apparently only works for smooth $Z$ and $Z'$, which is very restricted,
it is more than just a special case of (4.5) -- since it carries
certain nontrivial arithmetic information. \p

Now we begin the proof of Theorem B. Let $f\colon X \cdot\!\cdot\!\!\to X'$
be a birational map between two $n$ dimensional smooth projective 
varieties where only $X$ is assumed to be minimal. In the notation of
\S1, Theorem 1.4 says that $E\ge E'$. So we obtain canonical morphisms 
$H^i(E)\to H^i(E')$ induced from $E'\subset E$. Since $Z:=\phi(E)$ and 
$Z':=\phi'(E')$ are of codimension at least two, 
$H^{2n-2}(Z)=0=H^{2n-2}(Z')$. By comparing $(4.4)$ and $(4.4')$ 
via the surjective map $H^{2n-2}(E)\to H^{2n-2}(E')$, we obtain 
a canonical embedding:
$$
\phi^* H^{2n-2}(X,\Z)\subset \phi'^*H^{2n-2}(X',\Z). \leqno (4.7)
$$
which respects the Hodge structures. This induces an injective map
$$
\phi'_!\circ \phi^*: H^{2n-2}(X,\Z) \to H^{2n-2}(X',\Z), \leqno (4.8)
$$
which by the projection formula is easily seen to be independent of the
choice of $Y$, hence canonical. Poincar\'e duality then concludes 
the first statement of Theorem B.\par 

For the second statement, we may simply copy the above proof by replacing 
(4.4) with the similar formula for the Weil divisors. Q.E.D.\p
 
One can also interpret this result in terms of the Picard group
if the terminal varieties considered are assumed to be factorial or  
$\Q$-factorial.\p\p\p

\noindent
{\bf \S5 Further Remarks}\p

We conclude this paper with two historical remarks and two technical
remarks:\p

{\bf 5.1.} A version of Theorem 1.4, or rather the Corollary
1.10, was used before by Koll\'ar in his study of three dimensional 
flops. In fact, he proved that three dimensional birational 
$\Q$-factorial minimal models all share the same singularities, 
singular cohomologies and intersection cohomologies with pure Hodge 
structures (for deep reasons). See [Ko] for the details.\par

More recently, the author used a relative version of (1.10) to study 
degenerations of minimal projective threefolds [Wa; \S4] and obtained a 
negative answer to the so called ``filling-in problem'' 
in dimension three. Namely, there exist degenerating projective families 
of smooth threefolds which are $C^\infty$ trivial over the punctured disk, 
but can not be completed into smooth projective families.\p 

{\bf 5.2.} After Koll\'ar's result on threefolds, the problem 
on the equivalence of Betti numbers seemed to be ignored for a while
until recently when Batyrev treated the case of projective 
Calabi-Yau manifolds [Ba].\par

In the special case of projective hyper-K\"ahler 
manifolds, Theorem A has also been proved recently by Huybrechts [Hu] 
using quite different methods. In fact, he proves more -- these manifolds 
are all inseparable points in the moduli space (hence are diffeomorphic 
and share the same Hodge structures)!\par

This problem on general minimal models, to the best of the author's
knowledge, has not been studied before our paper. In our case, the homotopy 
types will generally be different. In fact, it is well known that for 
a single elementary 
transform of threefolds, although the singular cohomologies are 
canonically identified, the cup product must change. However, inspired
by Koll\'ar's result and Remark 3.7, we still expect that the (non-polarized) 
Hodge structures will turn out to be the same.\p

{\bf 5.3.} In order to generalize Theorem A to the singular case, 
our approach works equally well in the log-terminal case, with 
the only problem being that we need a good interpretation like 
Weil's formula (2.8) for the precise meaning of the weighted counting, 
which is the key to relate 
$p$-adic integrals to the Weil conjecture.\par
 
Since a suitable version of the Weil conjecture for
singular varieties has already been proved by Deligne in [BBD] in terms
of the intersection cohomologies introduced by Goresky and MacPherson
[GM], this problem is thus reduced to the calculation of local Lefschetz 
numbers. \par

More precisely, one needs to evaluate the $p$-adic integrals
over a singular point and to reconstruct the ``constructible complexes
of sheaves'' which it may correspond to. If luckily enough, it is the 
intersection cohomology complexes, then we may get our conclusion again 
via Deligne's theorem. A detailed discussion on this will be continued in a 
subsequent paper.\p

{\bf 5.4.} For Theorem B, it is likely that a similar argument
works for proving that terminal minimal models also minimize the second
intersection cohomology groups and that they all share the same pure Hodge 
structures. The important injectivity of $\phi^*: I\!H^i(X) \to I\!H^i(Y)$
needed to conclude (4.4) is now a consequence of the so called 
``decomposition theorem'' of projective morphisms. ([BBD] again!)\par

An interesting question arises: is the Picard number 
(or the second Betti number) of a non-minimal model always strictly 
bigger than the one attained by the minimal models?\par

Mazur raised the following question: can one extract the expected 
``minimal cohomology piece'' directly from any smooth model 
without refering to the minimal models?   
\p\p\p

\noindent
{\bf Acknowledgement}\p

During the preparation of this paper, I benefitted a lot from many 
interesting discussions with W.T. Gan, K. Zuo and B. Cui. Professor
C.-L. Chai pointed out several inaccuracies and provided useful comments
on the preliminary version. B. Hassett and B. Cheng carefully read the 
preliminary version and suggested improvements on the presentation.
I would like to take this opportunity to thank all of them.\par

I am also grateful to Professors R. Bott, H. Esnault, B. Gross, J. Harris, 
B. Mazur, W. Schmid, E. Viehweg and J. Yu for their interest in this work.
Last but not the least, I want to thank Professor S.T. Yau
for his steady encouragement.
\p\p\p

\noindent
{\bf References}\p

\item{[Ba]} V. Batyrev, {\it On the Betti numbers of birationally
isomorphic projective varieties with trivial canonical bundles},
alg-geom/9710020.
\item{[BBD]} A. Beilinson, J. Bernstein, P. Deligne, {\it Faiscuax
Pervers}, Ast\'erisque {\bf 100}, Soc.\ Math.\ de France 1983.
\item{[Ca]} W.S. Cassels, {\it Local Fields}, London Math. Soc. Student
 Text {\bf 3}, Cambridge Univ.\ Press 1986. 
\item{[D1]} P. Deligne, {\it La conjecture de Weil I}, IHES Publ.\ Math.\
{\bf 43} (1974), 273-307. {\it II}, {\bf 52} (1980), 313-428.
\item{[D2]} ---{}--- {\it et al.}, {\it S\`eminaire de G\`eom\`etrie
Alg\`ebrique du Bois-Marie} SGA $4{1\over 2}$, Lect.\ Notes in Math.\ 
{\bf 569}, Springer-Verlag, NY 1977. 
\item{[GM]} M. Goresky, R. MacPherson, {\it Intersection homology II},
Inv.\ Math.\ {\bf 71} (1983), 77-129.
\item{[HN]} G. Harder, M.S. Narasimhan, {\it On the cohomology groups
of moduli spaces of vector bundles over curves}, Math.\ Annln.\ {\bf 212}
(1975), 215-248.
\item{[Hi]} H. Hironaka, {\it Resolution of singularities of an algebraic
variety over a field of characteristic zero}, Ann.\ of Math.\ {\bf 79} 
(1964), 109-326.
\item{[Hu]} D. Huybrechts, {\it Compact hyperk\"ahler manifolds: basic 
results}, alg-geom/9705025.
\item{[Ig]} J.-I. Igusa, {\it Lecture on Forms of Higher Degree}, Tata
Inst.\ of Fund.\ Res., Bombay 1978.
\item{[KMM]} Y. Kawamata, K. Matsuda, K. Matsuki, {\it Introduction to
the minimal model program}, Adv.\ Stud.\ in Pure Math.\ {\bf 10} (1987), 
283-360.
\item{[Ko]} J. Koll\'ar, {\it Flops}, Nagoya Math.\ J. {\bf 113} (1989),
15-36.
\item{[Wa]} C.-L. Wang {\it On the incompleteness of the Weil-Petersson
metric along degenerations of Calabi-Yau manifolds}, Math.\ Res.\ Let.\ 
{\bf 4} (1997), 157-171.
\item{[W1]} A. Weil, {\it Numbers of solutions of equations in finite
fields}, Bull.\ AMS {\bf VI.55} (1949), 497-508.
\item{[W2]} ---{}---, {\it Ad\`ele and Algebraic Groups}, Prog.\ Math.\ 
{\bf 23}, Birkhauser, Boston 1982.
\p\p\p

Chin-Lung Wang\par
Harvard University, Department of Mathematics\par
1 Oxford Street 325, Cambridge, MA 02138\par
Email: dragon@math.harvard.edu\par 

\end